\documentstyle[amssymb]{amsart} 

\reversemarginpar

\reversemarginpar

\reversemarginpar
\newcommand{\baire}{\omega^{\textstyle \omega}}
\newcommand{\can}{2^{\textstyle \omega}}
\newcommand{\non}{\text{\normalshape\sf{non}}}
\newcommand{\df}{{\text{\normalshape\sf {def}}}}
\newcommand{\lh}{{\text{\normalshape\sf {lh}}}}
\newcommand{\fs}{2^{\textstyle <\!\omega}}
\newcommand{\conc}{{}^\frown\!}
\newcommand{\gb}{{\mathfrak b}}
\newcommand{\gd}{{\mathfrak d}}

\newcommand{\reals}{{\Bbb R}}

\newcommand{\rest}{{\mathord{\restriction}}}

\newcommand{\cov}{\text{\normalshape\sf{cov}}}
\newcommand{\unif}{\text{\normalshape\sf{non}}}
\newcommand{\Borel}{\text{\normalshape\sf{Borel}}}

\newcommand{\dom}{{\text{\normalshape\sf {dom}}}}

\newcommand{\nor}{|}

\newcommand{\QED}{\hspace{0.1in} \Box \vspace{0.1in}}
\newcommand{\F}{{\cal F}}

\newcommand{\N}{{\cal N}}
\newcommand{\M}{{\cal M}}
\newcommand{\V}{{\bold V}}

\renewcommand{\>}{\rangle}

\newtheorem{theorem}{Theorem}[section]
\newtheorem{lemma}[theorem]{Lemma}
\newtheorem{corollary}[theorem]{Corollary}
\newtheorem{claim}{Claim}[theorem]
\newtheorem{proposition}[theorem]{Proposition}

\theoremstyle{definition}
\newtheorem{definition}[theorem]{Definition}
\newtheorem{notation}[theorem]{Notation}

\newcommand{\lesdot}{\mathrel{\mathord{<}\!\!\raise 0.8 pt\hbox{$\scriptstyle\circ$}}}
\newcommand{\Proof}{\noindent{\sc Proof} \hspace{0.2in}}

\newcommand{\lft}[2]{\mathopen\ifcase#1{}\oo\or
                        \big#2\or\Big#2\else\oo\fi} 
\newcommand{\rgt}[2]{\mathclose\ifcase#1{}\oo\or
                        \big#2\or\Big#2\else\oo\fi}

\newcommand{\SM}{{\mathcal {SM}}}
\newcommand{\SN}{{\mathcal {SN}}}

\begin{document}
\title{After all, there are some inequalities which are provable in ZFC}
\author{Tomek Bartoszy\'{n}ski}
\thanks{First author partially supported by  NSF grant DMS 95-05375
  and Alexander von Humbold Foundation}
\address{Department of Mathematics\\
Boise State University\\
Boise, Idaho 83725 U.S.A.}
\email{tomek@@math.idbsu.edu, http://math.idbsu.edu/\~{}tomek/}
\author{Andrzej Ros{\l}anowski}
\thanks{Second author partially supported by KBN grant 2 PO3 A 01109}
\address{Department of Mathematics\\
Hebrew University\\
Jerusalem, Israel, and\\
Mathematical Institute\\
Wroclaw University\\
Wroclaw, Poland}
\email{roslanow@@sunrise.huji.ac.il, http://www.ma.huji.a.il/\~{}roslanow}
\author{Saharon Shelah}
\thanks{Third author partially supported by Basic Research Fundation
administered by Israel Academy of Sciences. Publication 616}
\address{Department of Mathematics\\
Hebrew University\\
Jerusalem, Israel}
\email{shelah@@sunrise.huji.ac.il, http://math.rutgers.edu/\~{}shelah/} 
\subjclass{03E35}

\begin{abstract}
We address ZFC inequalities between some cardinal invariants of the
continuum, which turned to be true in spite of strong expectations
given by \cite{RoSh:470}. 
\end{abstract}
\maketitle

\section{Introduction}
The present paper consists two independent sections which have two things in
common: both resulted in a  failure to fulfill old promises to build a
specific forcing notions and in both an important role is played by 
cardinal invariant $\kappa^*$. 

The first promise was stated in \cite{BRSh:490} and was related to cardinal
invariant $\cov^\star(\N)$. Let ${\bold B}$ denote the measure algebra adding
one random real.  

\begin{definition}
Let $\N_2$ be the ideal of measure zero subsets of $\reals\times\reals$ and
let $\Borel(\reals)$ be the collection of all Borel mappings from $\reals $
into $\reals$. Define 
\begin{multline*}
\cov^\star(\N)=\min\lft2\{\nor{\cal A}\nor:  {\cal A} \subseteq \N_2 \ \&\
(\forall f\in\Borel(\reals))(\forall B\in\Borel\setminus\N)(\exists H \in
{\cal A})\\ 
(\lft1\{x \in B: \<x,f(x)\> \in H\lft1\}\not\in \N)\rgt2\}
\end{multline*}
and 
\begin{multline*}
\unif^\star(\N)=\min\lft2\{\nor X\nor: X\subseteq \Borel(\reals) \ \& \ (\forall
H\in\N_2)(\forall B\in\Borel\setminus\N)(\exists f\in X)\\
(\lft1\{x\in B : \<x,f(x)\>\not\in H\rgt1\}\not\in \N)\rgt2\}.
\end{multline*}
\end{definition}

\begin{proposition}\label{prop1}
$\cov^\star(\N)=\cov(\N)^{\V^{\bold B}}$ and
$\unif^\star(\N)=\unif(\N)^{\V^{\bold B}}$.~$\QED$ 
\end{proposition}
It has been known that (see \cite{BRSh:490}, \cite{Kra83}, \cite{Paw86Sol} for
more details): 
\begin{enumerate}
\item $\max\{\cov(\N)^{\V},{\frak b}^{\V}\}\leq\cov(\N)^{\V^{\bold B}}\leq
\non(\M)$; 
\item it is consistent that $\cov(\N)^{\V^{\bold B}}>\max\{\cov(\N)^{\V},
{\frak b}^{\V}\}$;
\item it is consistent that $\cov(\N)^{\V^{\bold B}}> {\frak d}$.
\end{enumerate}
And in \cite[3.11]{BRSh:490} we promised that in \cite{RoSh:470} it would be
proved that 
\begin{itemize}
\item it is consistent that $\cov(\N)^{\V^{\bold B}}<\non(\M)$,
\end{itemize}
being sure that using the method of norms on possibilities we could construct a
forcing notion which:
\begin{description}
\item[a] is proper $\baire$--bounding,
\item[b] makes ground reals meager and
\item[c] does not add a $\bold B$--name for a random real over $\V^{\bold B}$.
\end{description}
However, when trying to fill up the details of the construction, we have
discovered that {\em there is no such forcing notion} and found new
inequalities provable in ZFC.

The second section deals with an inequality related to localizations of
subsets of $\omega$ by partitions of $\omega$. Several notions of localization
and related cardinal invariants were introduced in \cite{RoSh:501}. The one we
will refer to is the $R^\exists_0$--localization property.
\begin{definition}
\label{localizations}
Let $\V\subseteq\V^*$ be universes of Set Theory and let $k\in\omega$. 
\begin{enumerate}
\item We say that the extension $\V\subseteq\V^*$ has {\em
the $R^\exists_k$--localization property} if in $\V^*$:
\begin{quotation}
\noindent for every infinite co-infinite set $X\subseteq\omega$ there is a
partition $\langle K_n: n\in\omega\rangle\in V$ of $\omega$ such that
$\nor K_n\nor >k+1$ and 
\[(\exists^\infty n\in\omega)(\nor X\cap K_n\nor \leq k).\]
\end{quotation}
\item An infinite co-infinite set $X\subseteq\omega$, $X\in\V^*$ is said to be
{\em $(k,0)$--large over $\V$} if
\begin{quotation}
\noindent for every sequence $\langle K_n: n\in\omega\rangle\in \V$ of
disjoint $k$--element subsets of $\omega$ we have
\[(\forall^\infty n\in\omega)(K_n\cap X\neq\emptyset).\]
\end{quotation}
\end{enumerate}
\end{definition}
The following result has been shown in \cite[1.8]{RoSh:501}.

\begin{proposition}
\label{thereisno}
Let $\V\subseteq\V^*$ be models of ZFC, $m\geq 2$, $k\in\omega$. Then the
following conditions are equivalent:
\begin{enumerate}
\item there is no $(m,0)$-large set in $\V^*$ over $\V$,
\item there is no $(2,0)$-large set in $\V^*$ over $\V$,
\item $\V\subseteq\V^*$ has the $R^\exists_0$--localization property,
\item $\V\subseteq\V^*$ has the $R^\exists_k$--localization property.~$\QED$
\end{enumerate}
\end{proposition}
After noting that if $\V\cap\can$ is not meager in $\V^*$, $\V\subseteq\V^*$
then the extension $\V\subseteq\V^*$ has the $R^\exists_0$--localization
property we promised to give in \cite{RoSh:470} an example of a forcing notion
showing that the converse implication does not hold. In fact we wanted to
construct a forcing notion which:
\begin{description}
\item[a] is proper $\baire$--bounding,
\item[b] makes ground reals meager and
\item[c] has the $R^\exists_0$--localization property.
\end{description}
Once again, we have discovered that {\em there is no such forcing notion} and
we have established some new inequalities between relevant cardinal
invariants. 


\begin{notation}
We try to keep our notation standard and compatible with that of classical
textbooks on Set Theory (like Jech \cite{J} or Bartoszy{\'n}ski Judah
\cite{BJbook}) 
) 
\begin{enumerate}
\item Let $i,j<\omega$. The set of all integers $m$ satisfying $i\leq m<j$ is
denoted by $[i,j)$, etc.
\item For integers $k_i,\ldots,k_j$ ($i\leq j<\omega$), $\prod\limits_{\ell
=i}^j k_\ell$ is their Cartesian product interpreted as the collection of all
finite functions $\tau$ such that $\dom(\tau)=[i,j]$ and $(\forall\ell\in
\dom(\tau))(\tau(\ell)\in k_i)$. 

However, we will use the same notation for the cardinality of this set, hoping
that it will not cause too much confusion.
\item For two sequences $\eta,\nu$ we write $\nu\vartriangleleft\eta$ whenever
$\nu$ is a proper initial segment of $\eta$, and $\nu\trianglelefteq\eta$ when
either $\nu\vartriangleleft\eta$ or $\nu=\eta$. The length of a sequence
$\eta$ is denoted by $\lh(\eta)$.
\item The quantifiers $(\forall^\infty n)$ and $(\exists^\infty n)$ are
abbreviations for  
\[(\exists m\in\omega)(\forall n>m)\quad\mbox{ and }\quad(\forall m\in\omega)
(\exists n>m),\] 
respectively.
\item For $\omega$--sequences $\eta,\rho$ we write $\eta=^*\rho$ whenever
\[(\forall^\infty n\in\omega)(\eta(n)=\rho(n)).\]
\item The Cantor space $\can$ and the Baire space $\baire$ are the spaces of
all functions from $\omega$ to $2$, $\omega$, respectively, equipped with
natural (product) topology. 
\item In forcing arguments, {\em a stronger condition is the larger one}.
\end{enumerate}
\end{notation}

\section{Adding a random name for a random real}
As the failure in building the forcing notion we had in mind for 
\cite[3.11]{BRSh:490} directly results in some properties of extensions of
universes of ZFC, we will formulate the main result of the present section in
this language. Further we will draw several conlusions for cardinal
invariants. 

The result presented in \ref{addname} below  is of some interest {\em per se} if you
have in mind the following theorem (see \cite[3.1]{Paw86Sol}).

\begin{theorem}
1.\ \ (Krawczyk; \cite{Kra83})\ \ \ Suppose that $\V\subseteq\V^*$ are
universes of Set Theory such that $\V\cap\baire$ is bounded in
$\V^*\cap\baire$. Let $r$ be a random real over $\V^*$. Then there is in
$\V^*[r]$ a random real over $\V[r]$.\\
2.\ \ (Pawlikowski; \cite[3.2]{Paw86Sol})\ \ \ Suppose that $c$ is a Cohen
real over $\V$ and $r$ is a random real over $\V[c]$. Then, in $\V[c][r]$
there is no random real over $\V[r]$. 
\end{theorem}

\begin{definition}
\label{constructor}
Let $\Phi\in\baire$ be a strictly increasing function. {\em A
$\Phi$--constructor} is a sequence $\langle n_i,m_i,k_i:i<\omega\rangle$ of
integers defined inductively by:\quad $n_0=2$ and for $i\in\omega$
\[m_i=(\prod\limits_{j<i} m_j)\cdot 2^{3(n_i+i)},\quad k_i=(m_i
\cdot\prod\limits_{j<i} k_j)\cdot\Phi(m_i\cdot\prod\limits_{j<i} k_j),\quad
n_{i+1}=n_i(k_i+1).\]
[So $n_i<m_i<k_i<n_{i+1}$.]
\end{definition}

\begin{theorem}
\label{addname}
Suppose that $\V\subseteq\V^*$ are universes of Set Theory such that 
\begin{quotation}
\noindent if $r$ is a random real over $\V^*$

\noindent then in $\V^*[r]$ there is no random real over $\V[r]$.
\end{quotation}
Let $\Phi\in\baire\cap\V$ be strictly increasing and let $\langle n_i,m_i,k_i:
i<\omega\rangle$ be the $\Phi$--constructor. 
\begin{enumerate}
\item Assume $\V\cap\baire$ is dominating in $\V^*\cap\baire$. Then, in $\V^*$:

for every function $\eta\in\prod\limits_{i\in\omega}k_i$ there are sequences
$\langle X_\ell:\ell<\omega\rangle\in\V$ and $\langle
i_m:m<\omega\rangle\in\V$ such that
\begin{description}
\item[a] $(\forall \ell\in\omega)(X_\ell\subseteq k_\ell\ \&\ \nor X_\ell\nor =
m_\ell\
\cdot\prod\limits_{j<\ell} k_j)$,
\item[b] $(\forall m\in\omega)(\exists \ell\in [i_m,i_{m+1}))(\eta(\ell) \in
X_\ell)$. 
\end{description}
\item Assume $\V\cap\baire$ is unbounded in $\V^*\cap\baire$. Then, in $\V^*$:

for every function $\eta\in\prod\limits_{i\in\omega}k_i$ there is a sequence
$\langle X_\ell\!:\ell\!<\!\omega\rangle\in\V$ such that
\begin{description}
\item[a] $(\forall \ell\in\omega)(X_\ell\subseteq k_\ell\ \&\ \nor X_\ell\nor =
m_\ell\cdot\prod\limits_{j<\ell} k_j)$,
\item[b] $(\exists^\infty \ell\in \omega)(\eta(\ell)\in X_\ell)$. 
\end{description}
\end{enumerate}

\end{theorem}

\begin{theorem}
\label{addname1}
Suppose that $\V\subseteq\V^*$ are universes of Set Theory such that 
\begin{quotation}
\noindent if $r$ is a random real over $\V^*$

\noindent then  $\V[r] \cap 2^\omega$ has measure zero in  $\V^*[r]$.
\end{quotation}
Let $\Phi\in\baire\cap\V$ be strictly increasing and let $\langle n_i,m_i,k_i:
i<\omega\rangle$ be the $\Phi$--constructor. 
\begin{enumerate}
\item Assume $\V\cap\baire$ is bounded in $\V^*\cap\baire$. Then
 there are sequences
$\langle X_\ell:\ell<\omega\rangle\in\V^*$ and $\langle
i_m:m<\omega\rangle\in\V^*$ such that for every function
$\eta\in\prod\limits_{i\in\omega}k_i \cap \V$
\begin{description}
\item[a] $(\forall \ell\in\omega)(X_\ell\subseteq k_\ell\ \&\ \nor X_\ell\nor =
m_\ell\
\cdot\prod\limits_{j<\ell} k_j)$,
\item[b] $(\forall^\infty  m\in\omega)(\exists \ell\in [i_m,i_{m+1}))(\eta(\ell) \in
X_\ell)$. 
\end{description}
\item Assume $\V\cap\baire$ is unbounded in $\V^*\cap\baire$ (but not
  dominating). Then, there is a sequence
$\langle X_\ell\!:\ell\!<\!\omega\rangle\in\V^*$ such that

for every function $\eta\in\prod\limits_{i\in\omega}k_i \cap \V^*$ 
\begin{description}
\item[a] $(\forall \ell\in\omega)(X_\ell\subseteq k_\ell\ \&\ \nor X_\ell\nor =
m_\ell\cdot\prod\limits_{j<\ell} k_j)$,
\item[b] $(\exists^\infty \ell\in \omega)(\eta(\ell)\in X_\ell)$. 
\end{description}
\end{enumerate}

\end{theorem}

{\sc Proof of \ref{addname}} We will only prove \ref{addname}, the proof of \ref{addname1}
is obtained by dualization.

The main parts of the proofs of (1) and (2) are the same, the
difference comes only at the very end. So, for a while, we will not specify
which part of the theorem we are proving. We will present a construction which
itself is interesting, though it is very elementary. 

Let $\Phi\in\baire\cap\V$ be increasing and let $\langle n_i,m_i,k_i:
i<\omega\rangle$ be the $\Phi$--constructor. Letting $n_{-1}=0$, for each
$i\in\omega$ choose a sequence $\langle f^i_\ell: \ell<k_i\rangle$ of
functions such that
\begin{itemize}
\item $f^i_\ell: 2^{[n_i,n_{i+1})}\longrightarrow 2^{[n_{i-1},n_i)}$,
\item for every sequence $\langle\nu_\ell:\ell<k_i\rangle\subseteq 2^{[n_{i-1},
n_i)}$ we have
\[\nor \left\{\rho\in 2^{[n_i,n_{i+1})}: (\forall \ell<k_i)(f^i_\ell(\rho)=\nu_\ell)\right\}
\nor =\frac{2^{n_{i+1}-n_i}}{2^{(n_i-n_{i-1})k_i}}=2^{k_i\cdot n_{i-1}}.\]
\end{itemize}
For $i\leq j<\omega$ and $\eta\in\prod\limits_{r=i}^{j} k_r$ let
\[f^{i,j}_\eta: 2^{[n_i,n_{j+1})}\longrightarrow 2^{[n_{i-1},n_j)}:\ \ \rho
\mapsto \bigcup_{r=i}^j f^r_{\eta(r)}(\rho\rest [n_r,n_{r+1})).\]
The main point of our arguments will be done by the following combinatorial
observation (which should be clear if $S$ is thought of as a tree of
independent equally distributed random variables, but still it needs some
calculations). 

\begin{claim}
\label{probability}
Suppose that $0<i\leq j<\omega$ and $\emptyset\neq S\subseteq\prod\limits_{
r=i}^j k_r$ is such that for each $\eta\in S$ and $r\in [i,j]$:
\[\nor \{\tau(r): \tau\in S\ \&\ \tau\rest r=\eta\rest r\}\nor =m_r.\]
Then
\begin{multline*}
\left| \left\{\rho\in 2^{[n_i,n_{j+1})}\!: (\exists\sigma\!\in\! 2^{[n_{i-1},n_j)})
\left(\frac{7}{2^{n_j-n_{i-1}}}<\frac{\nor \{\tau\in S\!: f^{i,j}_\tau(\rho)=\sigma\}
\nor }{\nor S\nor }\right)\right\}\right| <\\
\frac{1}{2}\cdot 2^{n_{j+1}-n_i}.
\end{multline*}
\end{claim}

\noindent{\em Proof of the claim:}\hspace{0.15in} Fix $r\in [i,j]$ and $\tau^*
\in\prod\limits_{\ell\in [i,r)} k_\ell$ (so if $r=i$ then $\tau^*=\langle
\rangle$) such that there is $\tau\in S$ with $\tau^*\vartriangleleft\tau$. Let

\begin{multline*}
A^r_{\tau^*}\stackrel{\df}{=}\bigg\{\rho\in 2^{[n_r,n_{r+1})}:\ \mbox{ for some
}\sigma\in 2^{[n_{r-1},n_r)}\\
\left| \frac{\nor \{\tau(r):\tau\in S\ \&\ \tau^*\vartriangleleft\tau\ \&\ f^r_{\tau(r)}
(\rho)=\sigma\}\nor }{m_r}-\frac{1}{2^{n_r-n_{r-1}}}\right| \geq\frac{1}{2^{n_r-n_{r-1}}
\cdot 2^r}\bigg\}.
\end{multline*}
By Bernoulli's law of large numbers and by the definition of the $m_i$'s we
know that
\[\frac{\nor A^r_{\tau^*}\nor }{2^{n_{r+1}-n_r}}\leq 2^{n_r-n_{r-1}}\cdot \frac{1}
{4\cdot m_r\cdot (2^{-(n_r-n_{r-1}+r)})^2}= \frac{1}{4\cdot\prod\limits_{\ell<
r}m_\ell}\cdot 2^{-(3n_{r-1}+r)}.\]
Let 
\[A\stackrel{\df}{=}\left\{\rho\!\in\! 2^{[n_i,n_{j+1})}\!: (\exists r\!\in\!
[i,j])(\exists\tau^*\!\in\!\prod_{\ell\in [i,r)} k_\ell)\big((\exists\tau\!\in
\!S)(\tau^*\vartriangleleft\tau)\ \&\ \rho\rest [n_r,n_{r+1})\in A^r_{\tau^*}
\big)\right\}.\]
Note that
\begin{multline*}
\frac{\nor A\nor }{2^{n_{j+1}-n_i}}\leq\sum\limits_{r\in [i,j]}\sum\left\{\frac{\nor 
A^r_{\tau^*}\nor }{2^{n_{r+1}-n_r}}: \tau^*\in\prod\limits_{\ell\in [i,r)}
k_\ell\ \&\ (\exists\tau\in S)(\tau^*\vartriangleleft\tau)\right\}\leq\\
\sum\limits_{r\in [i,j]}\left(\frac{1}{4\cdot\prod\limits_{\ell<r}m_\ell}\cdot
2^{-(3n_{r-1}+r)}\cdot\prod\limits_{\ell<r}m_\ell\right)\leq\frac{1}{4} \sum
\limits_{r\in [i,j]} 2^{-r}<\frac{1}{2}.
\end{multline*}
Suppose now that $\rho\in 2^{[n_i,n_{j+1})}\setminus A$. Let $\sigma\in
2^{[n_{i-1},n_j)}$. We know that for each $r\in [i,j]$ and $\tau^*\in
\prod\limits_{\ell\in [i,r)}k_\ell$ such that $(\exists\tau\in S)(\tau^*
\vartriangleleft\tau)$ we have $\rho\rest [n_r,n_{r+1})\notin A^r_{\tau^*}$
and therefore 
\[\frac{\nor \{\tau(r): \tau\in S\ \&\ \tau^*\vartriangleleft\tau\ \&\
f^r_{\tau(r)}(\rho\rest [n_r,n_{r+1}))=\sigma\rest [n_{r-1},n_r)\}\nor }{m_r}<
(1+\frac{1}{2^r})\cdot \frac{1}{2^{n_r-n_{r-1}}}\]
(just look at the definition of the set $A^r_{\tau^*}$). Hence
\begin{multline*}
\frac{\nor \{\tau\in S: f^{i,j}_\tau(\rho)=\sigma\}\nor }{\nor S\nor }<\prod\limits_{r=
i}^j (1+\frac{1}{2^r})\frac{1}{2^{n_r-n_{r-1}}}=\frac{1}{2^{n_j-n_{i-1}}}
\cdot\prod\limits_{r=i}^j (1+\frac{1}{2^r})<\\
\frac{1}{2^{n_j-n_{i-1}}}\cdot e^{2^{1-i}} <\frac{7}{2^{n_j-n_{i-1}}}.
\end{multline*}
This finishes the proof of the claim.
\medskip

\noindent Now define a function 
\[F:\prod_{i\in\omega}k_i\times\can\longrightarrow\can:\ (\eta,\rho)\mapsto
\bigcup_{i\in\omega} f^i_{\eta(i)}(\rho\rest [n_i,n_{i+1})).\]
It should be clear that $F$ is well defined (look at the choice of the
$f^i_\ell$'s) and its definition (or rather its code) is in $\V$. The function
$F$ is continuous and we have the following claim.

\begin{claim}
\label{cl1}
If $\eta_0,\eta_1\in\prod_{i\in\omega}k_i$, $\rho_0,\rho_1\in\can$ and $\eta_0
=^*\eta_1$, $\rho_0=^*\rho_1$ then $F(\eta_0,\rho_0)=^*F(\eta_1,\rho_1)$. 
\end{claim}

\noindent{\em Proof of the claim:}\hspace{0.15in} Should be clear.
\medskip

\noindent Before we continue with the proof of the theorem let us introduce
some more notation. For a tree $T\subseteq\fs\times\fs$ and integers
$\ell,i<\omega$ we let $T_i\stackrel{\df}{=} T\cap (2^{n_{i+1}}\times
2^{n_i})$ and 
\begin{multline*}
T^{[\ell]}_i=\big\{(\nu_0,\nu_1)\in 2^{n_{i+1}}\times 2^{n_i}: \mbox{ if
$\ell<i$ then there are $(\nu_0',\nu_1')\in T_i$ such that}\\
\nu_0'\rest [n_{\ell+1},n_{i+1})=\nu_0\rest [n_{\ell+1},n_{i+1})\quad\mbox{
and }\quad\nu_1'\rest [n_\ell,n_i)=\nu_1\rest [n_\ell,n_i)\big\}.
\end{multline*}
If $\ell<i<\omega$ then we may treat members of $T^{[\ell]}_i$ as elements of
$2^{[n_{\ell+1},n_{i+1})}\times 2^{[n_\ell,n_i)}$ (as only this part carries
any information). Thus if $\rho_0\in 2^{[n_{\ell+1},n_{i+1})}$, $\rho_1\in
2^{[n_\ell,n_i)}$ then $(\rho_0,\rho_1)\in T^{[\ell]}_i$ means that there is
$(\nu_0,\nu_1)\in T^{[\ell]}_i$ such that $\nu_0\rest [n_{\ell+1},n_{i+1})=
\rho_0$, $\nu_1\rest [n_\ell,n_i)=\rho_1$. 

\begin{claim}
\label{cl2}
Suppose that $\eta\in\prod\limits_{i\in\omega}k_i\cap\V^*$. Then there is a
tree $T\subseteq\fs\times\fs$, $T\in \V$ such that
\begin{enumerate}
\item[(i)] $\mu^2([T])>0$

(where $[T]$ is the set of all infinite branches through $T$, 
\[[T]=\{(\rho,\sigma)\in\can\times\can: (\forall n\in\omega)((\rho\rest
n,\sigma\rest n)\in T)\},\]
and $\mu^2$ stands for the Lebesgue measure on the plane $\can\times\can$),
\item[(ii)] for each $\ell<\omega$
\[\mu\big(\{\rho\in\can: (\forall i\in\omega)((\rho\rest n_{i+1},F(\eta,\rho)
\rest n_i)\in T^{[\ell]}_i)\}\big)=0.\]
\end{enumerate}
\end{claim}

\noindent{\em Proof of the claim:}\hspace{0.15in} Let $r$ be a random real
over $\V^*$. By the assumptions of the theorem we know that $F(\eta,r)$ is not
a random real over $\V[r]$. Every Borel null subset of $\can$ from $\V[r]$ is
the section at $r$ of a Borel null subset of $\can\times\can$ from $\V$. 
Consequently we find a Borel null set $B\subseteq\can\times\can$ coded in $\V$
and such that $(r,F(\eta,r))\in B$. We may additionally require that $B$ is
invariant under rational translations, i.e. that 
\[(\rho_0,\rho_1)\in B\ \&\ \rho_0=^*\rho_0'\ \&\ \rho_1=^*\rho_1'\quad
\Rightarrow\quad (\rho_0',\rho_1')\in B.\]
In $\V$ take a closed subset of $\can\times\can$ of positive measure disjoint
from $B$. This gives a tree $T\in\V$, $T\subseteq\fs\times\fs$ such that
$\mu^2([T])>0$ and 
\begin{enumerate}
\item[$(\oplus)$] \qquad $(\forall\ell\in\omega)(\exists i\in\omega)((r\rest
n_{i+1},F(\eta,r)\rest n_i)\notin T^{[\ell]}_i)$. 
\end{enumerate}
We want to argue that this $T$ is as required and for this we need to check
the demand {\bf (ii)}. Let $\ell<\omega$. Look at the set 
\[Y\stackrel{\df}{=}\{\rho\in\can: (\forall i\in\omega)((\rho\rest n_{i+1},
F(\eta,\rho)\rest n_i)\in T^{[\ell]}_i)\}.\]
It is a closed subset of $\can$ coded in $\V^*$. Assume that $\mu(Y)>0$. Then
some finite modification $r^*$ of the random real $r$ is in $Y$. By \ref{cl1}
we know that $F(\eta,r)=^* F(\eta,r^*)$. Take $\ell_0>\ell$ so large that
\[F(\eta,r)\rest [n_{\ell_0},\omega)=F(\eta,r^*)\rest [n_{\ell_0},\omega)\quad
\mbox{ and }\quad r\rest [n_{\ell_0},\omega)=r^*\rest [n_{\ell_0},\omega).\]
Now note that 
\begin{multline*}
r^*\in Y\quad\Rightarrow\quad (\forall i\in\omega)((r^*\rest n_{i+1},F(\eta,
r^*)\rest n_i)\in T^{[\ell]}_i\quad\Rightarrow\\
\Rightarrow\quad (\forall i\in\omega)((r\rest n_{i+1},F(\eta,r)\rest n_i)\in
T^{[\ell_0]}_i
\end{multline*}
and the last contradicts $(\oplus)$ above, finishing the claim.

\begin{claim}
\label{cl3}
Suppose that $T\subseteq\fs\times\fs$ is a tree, $1<i<j<\omega$ and
\begin{enumerate}
\item[$(\otimes^i_j)$] \qquad $\frac{63}{64}\leq\frac{\nor T^{[i-1]}_j\nor }{2^{n_{j
+1}+n_j}}$.
\end{enumerate}
Let
\[W=\left\{\tau\in\prod_{\ell=i}^j k_\ell: \frac{\nor \{\rho\in 2^{[n_i,
n_{j+1})}: (\rho,f^{i,j}_\tau(\rho))\in T^{[i-1]}_j\}\nor }{2^{n_{j+1}-n_i}}<
\frac{1}{64}\right\}.\]
Then there are sets $X_i\subseteq k_i,\ X_{i+1}\subseteq k_{i+1},\ldots,
X_j\subseteq k_j$ such that
\begin{enumerate}
\item[$(\alpha)$] \quad $\nor X_\ell\nor \leq m_\ell\cdot \prod\limits_{r<\ell}k_r$
for each $\ell=i,\ldots,j$ and
\item[$(\beta)$]  \quad $(\forall\tau\in W)(\exists \ell\in [i,j])(\tau(\ell)
\in X_\ell)$.
\end{enumerate}
\end{claim}

\noindent{\em Proof of the claim:}\hspace{0.15in} Assume not. Then we may find
a set $S\subseteq \prod\limits_{\ell=i}^j k_\ell$ such that $S\subseteq W$ and
for every $\tau_0\in S$ and every $\ell\in [i,j]$
\[\nor \{\tau(\ell): \tau\in S\ \&\ \tau\rest \ell=\tau_0\rest \ell\}\nor =m_\ell.\]
How? For $\ell\in [i,j]$ let $W^\ell=\{\tau\rest\ell: \tau\in W\}$ (so $W^i
=\{\langle\rangle\}$). Now we choose inductively sets $X_\ell\subseteq k_\ell$
and $Y_\ell\subseteq W^\ell$ for $\ell=j,\ldots,i$. First we let 
\[Y_j=\big\{\nu\in W^j: \nor \{\tau(j):\nu\vartriangleleft\tau\in W\}\nor <m_j\big
\},\quad X_j=\bigcup_{\nu\in Y_j}\{\tau(j):\nu\vartriangleleft\tau\in W\}.\]
By its definition we have $\nor X_j\nor \leq m_j\cdot \nor Y_j\nor \leq m_j\cdot\prod
\limits_{r<j}k_r$. Suppose that $i\leq \ell<j$ and we have defined
$Y_{\ell+1} \subseteq W^{\ell+1}$ already. Let 
\[Y_\ell=\big\{\nu\in W^\ell: \nor \{\tau(\ell):\nu\vartriangleleft\tau\in
W^{\ell+1}\setminus Y_{\ell+1}\}\nor <m_\ell\big\},\qquad \mbox{ and}\]
\[X_\ell=\bigcup_{\nu\in Y_\ell}\{\tau(\ell):\nu\vartriangleleft\tau\in
W^{\ell+1}\setminus Y_{\ell+1}\}.\] 
Note that $\nor X_\ell\nor \leq m_\ell\cdot\nor Y_\ell\nor \leq m_\ell\cdot\prod\limits_{
r<\ell}k_r$.\\
Now look at the sets $X_i,\ldots,X_j$. By our assumption we know that there is
$\tau_0\in W$ such that $(\forall \ell\in [i,j])(\tau_0(\ell)\notin
X_\ell)$. This implies that $\langle\rangle\notin Y_i$. [Why? If
$\langle\rangle\in Y_i$ then, as $\tau_0(i)\notin X_i$, we have $\langle\tau_0
(i)\rangle\in Y_{i+1}$. Suppose have already shown that $\langle \tau_0(i),
\ldots,\tau_0(\ell)\rangle\in Y_{\ell+1}$, $i\leq\ell<j-1$. Since $\tau_0(\ell
+1)\notin X_{\ell+1}$ we conclude $\langle\tau_0(i),\ldots,\tau_0(\ell),
\tau_0(\ell+1)\rangle\in Y_{\ell+2}$. Thus, by induction, $\langle\tau_0(i),
\ldots,\tau_0(j-1)\rangle\in Y_j$ and $\tau_0(j)\in X_j$, a contradiction.]\\
Now we define the set $S\subseteq W$. We do this choosing inductively a finite
tree $S^*\subseteq\bigcup\limits_{\ell=i}^j \prod\limits_{r=i}^\ell k_r$ in
which maximal nodes will be elements of $W$. First we declare that $\langle
\rangle\in S^*$ and since $\langle\rangle\notin Y_i$ we may choose a set
$S^{\langle\rangle}_i\subseteq\{\tau(i): \tau\in W^{i+1}\setminus Y_{i+1}\}$ of
size $m_i$. We declare that $\{\langle z\rangle: z\in S^{\langle\rangle}_i\}
\subseteq S^*$. Note that $\langle z\rangle\in W^{i+1}\setminus Y_{i+1}$ for
$z\in S^{\langle\rangle}_i$. Suppose that we have decided that a sequence
$\nu\in \prod\limits_{r=i}^\ell k_r$ is in $S$, $i\leq \ell<j-1$ and we know
that $\nu\in W^{\ell+1}\setminus Y_{\ell+1}$. By the definition of
$Y_{\ell+1}$ we may choose a set $S^\nu_{\ell+1}\subseteq\{\tau(\ell+1):
\nu\vartriangleleft\tau\in W^{\ell+2}\setminus Y_{\ell+2}\}$ of size
$m_{\ell+1}$. We declare that all the sequences $\nu\conc\langle z\rangle$ for
$z\in S^{\nu}_{\ell+1}$ are in $S^*$. Note that we are sure that $\nu\conc
\langle z\rangle\in W^{\ell+2}\setminus Y_{\ell+2}$ (for $z\in S^{\nu}_{\ell+
1}$). Finally, having decided that a sequence $\nu\in W^j\setminus Y_j$ is in
$S^*$ we choose a set $S^\nu_j\subseteq\{\tau(j):\nu\vartriangleleft\tau\in
W\}$ of size $m_j$ and we declare $\nu\conc\langle z\rangle\in S^*$ for $z\in
S^\nu_j$. Immediately by the construction of $S^*$ we see that the set $S=S^*
\cap\prod\limits_{\ell=i}^j k_\ell$ is as required. 

\noindent Define:
\[u_0\stackrel{\df}{=}\left\{\rho\in 2^{[n_i,n_{j+1})}:\frac{\nor \{\tau\in S:
(\rho, f^{i,j}_\tau(\rho))\in T^{[i-1]}_j\}\nor }{\nor S\nor }\geq \frac{1}{8}\right\},\]
\[u_1\stackrel{\df}{=}\left\{\rho\in 2^{[n_i,n_{j+1})}:\frac{\nor \{\sigma\in
2^{[n_{i-1}, n_j)}: (\rho,\sigma)\in T^{[i-1]}_j\}\nor }{2^{n_j-n_{i-1}}}\leq
\frac{7}{8}\right\},\]
\[u_2\stackrel{\df}{=}\left\{\rho\in 2^{[n_i,n_{j+1})}\!: (\exists\sigma\!\in\!
2^{[n_{i-1},n_j)})\left(\frac{7}{2^{n_j-n_{i-1}}}<\frac{\nor \{\tau\in S\!: 
f^{i,j}_\tau(\rho)=\sigma\}\nor }{\nor S\nor }\right)\right\}.\]
Since $S\subseteq W$, by Fubini theorem, we have that 
\[\frac{\nor u_0\nor }{2^{n_{j+1}-n_i}}<\frac{1}{8}.\]
Now look at the assumption $(\otimes^i_j)$ on $T$: it implies that, by Fubini
theorem once again, 
\[\frac{\nor u_1\nor }{2^{n_{j+1}-n_i}}\leq\frac{1}{8}.\]
Finally, by \ref{probability}, we know that
\[\frac{\nor u_2\nor }{2^{n_{j+1}-n_i}}\leq\frac{1}{2}.\]
Consequently we find a sequence $\rho\in 2^{[n_i,n_{j+1})}\setminus (u_0\cup
u_1\cup u_2)$. Since $\rho\notin u_0\cup u_1$ we know that in the sequence 
\[\langle f^{i,j}_\tau(\rho):\tau\in S\rangle\]
less than $\frac{1}{8}\cdot 2^{n_j-n_{i-1}}$ many values (from  $2^{[n_{i-1},
n_j)}$) appear more than $\frac{7}{8}\cdot\nor S\nor $ times. This implies that
there is one value $\sigma\in 2^{[n_{i-1}, n_j)}$ which appears in this
sequence more than $\frac{7}{2^{n_j-n_{i-1}}}\cdot\nor S\nor $ times and therefore
$\rho\in u_2$, a contradiction finishing the proof of the claim.
\medskip

Now we may prove the theorem.

\noindent (1)\ \ \ Assume that $\V\cap\baire$ is dominating in $\V^*\cap
\baire$. Let $\langle n_i,m_i,k_i: i<\omega\rangle$ be the $\Phi$--constructor
and let $F:\prod\limits_{i\in\omega}k_i\times\can\longrightarrow\can$ be as
defined above. Suppose $\eta\in\prod\limits_{i\in\omega}k_i$.\\
By Claim \ref{cl2} we find a tree $T\subseteq\fs\times\fs$ from $\V$
satisfying the demands (i) and (ii) of \ref{cl2}. Let $\varphi\in\baire\cap
\V^*$ be such that for each $i\in\omega$ 
\[i<\varphi(i)\quad\mbox{and}\quad\frac{\nor \{\rho\in 2^{[n_i,n_{\varphi(i)+
1})}: (\rho,f^{i,\varphi(i)}_{\eta\rest [i,\varphi(i)]}(\rho))\in T^{[i-1]
}_{\varphi(i)}\}\nor }{2^{n_{\varphi(i)+1}-n_i}}<\frac{1}{64}.\]
Since $\V\cap\baire$ is dominating in $\V^*\cap\baire$ we find an increasing
sequence of integers $\langle i_m: m\in\omega\rangle\in\V$ such that
\begin{enumerate}
\item[$(\otimes)$] \qquad $\frac{63}{64}\leq\frac{\nor T^{[i_0-1]}_j\nor }{2^{n_{j
+1}+n_j}}$\quad for each $i_0<j<\omega$,
\item[$(\otimes^+)$] for each $m\in\omega$
\[\frac{\nor \{\rho\in 2^{[n_{i_m},n_{i_{m+1}})}: (\rho,f^{i_m,i_{m+1}-1}_{\eta
\rest [i_m,i_{m+1})}(\rho))\in T^{[i_m-1]}_{i_{m+1}-1}\}\nor }{2^{n_{i_{m+1}}-
n_{i_m}}}<\frac{1}{64}.\]
\end{enumerate}
[Note that to get $(\otimes^+)$ it is enough to require $\varphi(i_m)<i_{m+1}$
for each $m\in\omega$, what is easy to get as $\V\cap\baire$ is dominating.]\\
Now we construct, in $\V$, a sequence $\langle X_\ell: \ell<\omega\rangle$.\\
Fix $m\in\omega$ for a moment. Note that $(\otimes)$ implies $(\otimes^{i_m}_{
i_{m+1}-1})$ of \ref{cl3}. Let
\[W_m=\left\{\tau\in\prod_{\ell=i_m}^{i_{m+1}-1} k_\ell: \frac{\nor \{\rho\in 2^{
[n_{i_m},n_{i_{m+1}})}: (\rho,f^{i_m,i_{m+1}-1}_\tau(\rho))\in T^{[i_m-1]}_{
i_{m+1}-1}\}\nor }{2^{n_{i_{m+1}}-n_{i_m}}}<\frac{1}{64}\right\}.\]
It follows from \ref{cl3} that there are sets $X_{i_m}\subseteq k_{i_m},
\ldots,X_{i_{m+1}-1}\subseteq k_{i_{m+1}-1}$ such that 
\begin{enumerate}
\item[$(\alpha)$] \quad $\nor X_\ell\nor \leq m_\ell\cdot \prod\limits_{r<\ell}k_r$
\item[$(\beta)$]  \quad $(\forall\tau\in W_m)(\exists \ell\in [i_m,i_{m+1}))
(\tau(\ell)\in X_\ell)$.
\end{enumerate}
But now we easily finish notifying that $(\otimes^+)$ implies that 
\[(\forall m\in\omega)(\eta\rest [i_m,i_{m+1})\in W_m).\]
\smallskip

\noindent (2)\ \ \ We repeat the arguments from the first case, but now we
cannot require $(\otimes^+)$. Still, as $\V\cap\baire$ is unbounded in $\V^*
\cap\baire$ we may demand that the sequence $\langle i_m:m\in\omega\rangle\in
\V$ satisfies $(\otimes)$ and
\begin{enumerate}
\item[$(\otimes^-)$] for infinitely many $m\in\omega$
\[\frac{\nor \{\rho\in 2^{[n_{i_m},n_{i_{m+1}})}: (\rho,f^{i_m,i_{m+1}-1}_{\eta
\rest [i_m,i_{m+1})}(\rho))\in T^{[i_m-1]}_{i_{m+1}-1}\}\nor }{2^{n_{i_{m+1}}-
n_{i_m}}}<\frac{1}{64}.\]
\end{enumerate}
Then, defining $W_m$ as above, we will have
\[(\exists^\infty m\in\omega)(\eta\rest [i_m,i_{m+1})\in W_m),\]
and this is enough to get the conclusion of (2).~$\QED$

\begin{corollary}
\label{coraddname}
Suppose that $\V\subseteq\V^*$ are universes of Set Theory such that 
\begin{quotation}
\noindent if $r$ is a random real over $\V^*$

\noindent then in $\V^*[r]$ there is no random real over $\V[r]$
\end{quotation}
and $\V\cap\baire$ is unbounded in $\V^*\cap\baire$. Let $H\in\baire\cap\V$ be
an increasing function. Then:
\[\V^*\models (\forall f\in\prod\limits_{\ell\in\omega}H(\ell))(\exists g\in
\prod\limits_{\ell\in\omega} H(\ell)\cap\V)(\exists^\infty \ell\in \omega)(g(
\ell)=f(\ell)).\]
\end{corollary}
\Proof Define inductively a sequence $\langle n_i,m_i,x_i, y_i,k_i:
i\in\omega\rangle\in\V$: 
\begin{multline*}
n_0=2,\quad m_0=64,\quad x_0=64,\quad y_0=64+\prod_{\ell<64}H(\ell),\quad
k_0=64\cdot y_0,\\ 
n_{i+1}=n_i\cdot (k_i+1),\quad m_{i+1}=\big(\prod_{j\leq i}m_j\big)\cdot
2^{3(n_{i+1}+i+1)},\quad x_{i+1}=x_i+m_{i+1}\cdot\prod_{j\leq i}k_j,\\ 
y_{i+1}=y_i+x_{i+1}+\prod_{\ell\in [x_i,x_{i+1})} H(\ell),\quad k_{i+1}=
y_{i+1}\cdot\big(m_{i+1}\cdot\prod_{j\leq i}k_j\big).
\end{multline*}
Note that $y_{i+1}-y_i>x_{i+1}>m_{i+1}\cdot\prod\limits_{j\leq i}k_j -m_i\cdot
\prod\limits_{j<i}k_j$. Consequently we may choose a strictly increasing
function $\Phi\in\baire\cap\V$ such that $(\forall i\in\omega)(\Phi(m_i\cdot
\prod\limits_{j<i}k_j)=y_i)$. Now look at the definition of the sequence
$\langle n_i,m_i,k_i: i\in\omega\rangle$ -- clearly it is the
$\Phi$--constructor.\\
For $i\in\omega$ we have $\prod\limits_{\ell\in [x_{i-1},x_i)} H(\ell)\leq
k_i$ (we let $x_{-1}=0$ here). So we may take a one--to--one function
$\pi_i: \prod\limits_{\ell\in [x_{i-1},x_i)} H(\ell)\longrightarrow k_i$.\\
Now suppose $f\in\prod\limits_{\ell\in\omega}H(\ell)\cap\V^*$. Define $\eta\in
\prod\limits_{i\in\omega} k_i\cap\V^*$ by
\[(\forall i\in\omega)(\eta(i)=\pi_i(f\rest [x_{i-1},x_i))).\]
By \ref{addname}(2) we find a sequence $\langle X_\ell: \ell\in\omega\rangle\in
\V$ satisfying \ref{addname}(2)(a),(b) (for our $\eta$). Using the sequence
$\langle X_\ell: \ell\in\omega\rangle$ (and working in $\V$) we define a
function $g\in\prod\limits_{r\in\omega} H(r)\cap\V$. Fix $\ell\in\omega$ and
look at the set  
\[Y_\ell\stackrel{\df}{=}\left\{\tau\in\prod_{r\in [x_{\ell-1},x_\ell)} H(r):
\pi_\ell(\tau)\in X_\ell\right\}.\]
Since $\nor Y_\ell\nor \leq m_\ell\cdot\prod\limits_{j<\ell}k_j=x_\ell-x_{\ell-1}$,
we find $\sigma_\ell\in\prod\limits_{r\in [x_{\ell-1},x_\ell)} H(r)$ such that
\[(\forall\tau\in Y_\ell)(\exists r\in [x_{\ell-1},x_\ell))(\sigma_\ell(r)=
\tau(r)).\]
Next let $g\in\prod\limits_{r\in\omega}H(r)\cap\V$ be such that $g\rest
[x_{\ell-1},x_\ell)=\sigma_\ell$ (for $\ell\in\omega$). We finish noting that
if $\eta(\ell)\in X_\ell$ then $f\rest [x_{\ell-1},x_\ell)\in Y_\ell$ and
therefore for some $r\in [x_{\ell-1},x_\ell)$ we have
$g(r)=f(r)$. ~$\QED$

\section{$\cov^*(\N)$ and other cardinal invariants}
Results of the previous section allow us to compare $\cov^*(\N)$ to
other cardinal invariants. 

We will need several definitions.
Let $f, g \in \omega^\omega $ be two nondecreasing functions such that
$0<g(n)<f(n)$ for every $n$.
Let $S_{f,g}=\prod_n [f(n)]^{g(n)}$ and $S^*_{f,g}=\prod_n
[f(n)]^{g(n)}\times [\omega]^\omega $.
Define relations $R^\forall_{f,g}, R^\exists_{f,g}$ as 
$$\eta R^\exists_{f,g} S \iff \exists^\infty n \ \eta(n)\in S(n)$$
$$\eta R^\forall_{f,g} S \iff \forall^\infty n \ \eta(n)\in S(n)$$
for $\eta \in \prod_n f(n)$ and $S \in S_{f,g}$.
In case when $g(n)=1$ for all $n$ we will drop subscript $g$ and
define
$$\eta_0 R^\exists_{f} \eta_1 \iff \exists^\infty n \
\eta_0(n)=\eta_1(n)$$
for $\eta_0, \eta_1 \in S_f$.
The dual relation $R^\forall_f$ is not very interesting, so we
consider the following weaker relations $R^{**}_{f,g}$ and $R^{**}_f$
defined as
$$\eta R^{**}_{f,g} (S,K) \iff \forall^\infty n \ \exists m \in
[k_n,k_{n+1}) \ \eta(m) \in S(m),$$
for $\eta\in S_f$, $S\in S_{f,g}$ and $K=\{k_0<k_1<\dots\}\in
[\omega]^\omega $.
Finally define for a relation $R \subseteq A\times B$,
$$\gb(R)=\min\{|X|: X \subseteq A \ \&\ \forall y \in B \ \exists x \in
X \ \neg xRy\}$$
$$\gd(R)=\min\{|Y|: Y \subseteq B \ \&\ \forall x \in A \ \exists y \in
Y \  xRy\}.$$

For various independence results and techniques connected with these
invariants see \cite{RoSh:470}.

Using this terminology we can express the results of the previous
section as follows.
\begin{theorem} There  are $f,g \in \omega^\omega $ such that 
  $\cov^*(\N)\geq \gd(R^\exists_{f,g})$.
If $\cov^*(\N)\geq \gd$ then $\cov^*(\N)\geq \gd(R^{**}_{f,g})$.

Similarly,
$\unif^*(\N)\leq \gb(R^\exists_{f,g})$, and if
$\unif^*(\N)\leq \gb$ then $\unif^*(\N)\leq \gb(R^{**}_{f,g})$.
\end{theorem}
\Proof
This is a simple reformulation of Theorem \ref{addname}.
Fix an increasing function $\Phi \in \omega^\omega$.
Let $M$ be a model of size $\cov^*(\N)$ containing a witness for
$\cov^*(\N)$, and containing $\Phi$.
Since $\cov^*(\N) \geq \gb$ we can assume that $M \cap \omega^\omega $
is an unbounded family.
Let $\{n_i,m_i,k_i: i \in \omega\}  \in M$ be a $\Phi$-constructor.
Define
$f(n)=k_n$ and $g(n)=m_n \prod_{i<n} f(i)$.
By \ref{addname}, 
$$ \forall \eta \in S_f \ \exists S \in S_{f,g} \cap M \
\exists^\infty n \ \eta(n) \in S(n).$$
Thus $\gd(R^\exists_{f,g}) \leq |M| = \cov^*(\N)$.
Remaining parts of the theorem are proved in the same way by using \ref{addname1}.
It is note very hard to see that by simple diagonalization we can show
that for many triples $(h,f,g)$ we have
$\gb(R^\exists_h) \leq \gb(R^\exists_{f,g})$ and $\gd(R^\exists_h)
\geq \gd(R^\exists_{f,g})$ 

\begin{definition}
  
Let
$$\kappa^*=\sup\left\{\gd(R^\exists_f) :f\in (\omega\setminus\{0\})^{\textstyle
\omega}\right\} \quad \text{and} \quad  \lambda^*=\inf\left\{\gb(R^\exists_f) :f\in (\omega\setminus\{0\})^{\textstyle
\omega}\right\}.$$
\end{definition}

\begin{theorem}
  $\cov^*(\N) \geq \kappa^*$ and $\unif^*(\N) \leq \lambda^*$.
\end{theorem}
\Proof
Let $f\in (\omega\setminus\{0\})^{\textstyle \omega}$. We may assume
that $f$ is strictly increasing. Take a family ${\cal A}\subseteq\N_2$
realizing the minimal cardinality in the definition of $\cov^*(\N)$ and take
an unbounded family $\F\subseteq\baire$ of size ${\frak b}$ (remember ${\frak
b}\leq\cov^*(\N)$). Let $N\prec({\cal H}(\chi),\in,<^*_\chi)$ be an elementary
submodel of size $\cov^*(\N)$ containing all members of ${\cal A}$ and $\F$
and such that $f\in N$. Now apply \ref{coraddname} to $N\subseteq\V$. Note
that if $r$ is a random real over $\V$ then in $\V[r]$ there is no random real
over $N[r]$ (as ${\cal A}\subseteq N$). Morever $N\cap\baire$ is unbounded in
$\V\cap\baire$ (as $\F\subseteq N$). Consequently (in $\V$) we have
\[(\forall h\in\prod_{n\in\omega}f(n))(\exists g\in \prod_{n\in\omega}f(n)\cap
N)(\exists^\infty n\in\omega)(g(n)=h(n)),\] 
showing that $\gd(R^\exists_f)\leq \nor N\nor =\cov^*(\N)$.~$\QED$

\begin{definition}
  Suppose that $ X \subseteq 2^\omega $. 
  \begin{enumerate}
  \item $X \in \SN$ (strong measure zero) if for every meager set $F
    \subseteq 2^\omega $, $X+F \neq 2^\omega $,
  \item $X \in \SM$ (strongly meager) if for every null set $H
    \subseteq 2^\omega $, $X+H \neq 2^\omega $,
  \end{enumerate}
\end{definition}

\begin{lemma}
  $ \lambda^*=\non(\SN)$ and $ \kappa^* \geq \non(\SM)$.
\end{lemma}
\Proof
The first equality was proved by Miller (see \cite{Mil81Som} or
\cite{BJbook}, 8.1.14).

Suppose that a family $\F \subseteq \prod_{n\in\omega} f(n)$
exemplifies $\gd(R^\exists_f)$.
Work in the space $X=\prod_{n\in\omega} f(n)$ (for sufficiently big $f$)
equipped with the standard product measure.
Consider the set $G=\{x \in X: \exists^\infty n
\ x(n)=0\}$. It is easy to see that $G$ is a null set and $\F+G=X$. 
Thus $\F \not\in \SM$ in $X$ (which easily translates to $2^\omega
$).~$\QED$

\begin{corollary}\label{corcovstar}
  $\cov^*(\N) \geq \max\{ {\mathfrak b}, \non(\SM)\}$ and
$\non^*(\N) \leq \min\{ {\mathfrak d}, \non(\SN)\}$.~$\QED$
\end{corollary}

\begin{lemma}
  If $ \cov^*(\N) \geq \gd$ then $\cov^*(\N)=\unif(\M)$.
  If $ \unif^*(\N) \leq \gb$ then $\unif^*(\N)=\cov(\M)$.
\end{lemma}
\Proof
We will prove only the first assertion. The other one is proved by the
dual argument.

It is well known (see \cite{BJbook}, 2.4.7, 2.4.1) that
$$\unif(\M)=\min\{|F|: F \subseteq \omega^\omega \ \& \ \forall g  \in
\omega^\omega \ \exists f \in F \ \exists^\infty n \ f(n)=g(n)\}$$
and 
$$\cov(\M)=\min\{|F|: F \subseteq \omega^\omega \ \& \ \forall g  \in
\omega^\omega \ \exists f \in F \ \forall^\infty  n \ f(n)\neq g(n)\}$$



Let $F \subseteq \omega^\omega $ be a dominating family of size
$\gd$. For each $f \in F$ choose a witness $X_f \subseteq S_f$ of size
$\gd(R^\exists_f)$.
Let $X=\bigcup_{f\in F} X_f$.
It is clear that $|X|=\max\{\gd,\kappa^*\} \leq \cov^*(\N)$ and
$$ \forall g \in \omega^\omega \ \exists f \in F\ \exists x_f \in X_f\
\exists^\infty n \ g(n)=x_f(n).$$
Thus, $\unif(\M) \leq \cov^*(\N)$. 
To see that $\cov^*(\N) \leq
\unif(\M)$ in we need the following lemma:
\begin{lemma}
$\cov^\star(\N) \leq \unif(\M)$ and $\unif^*(\N) \geq \cov(\M)$.  
\end{lemma}
\Proof
We have the following
$\cov^\star(\N)=\cov(\N)^{\V^{\bold B}} \leq \unif(\M)^{\V^{\bold
    B}}=\unif(\M)$. 
The first equality is by \ref{prop1}, the second is well known, and
for the third one see \cite{BJbook} or \cite{BRSh:490}.~$\QED$

\begin{corollary}
  There is no proper forcing notion $ {\mathcal P} $ such that 
\begin{enumerate}
\item is proper $\baire$--bounding,
\item makes ground reals meager and
\item does not add a $\bold B$--name for a random real over $\V^{\bold B}$.
\end{enumerate}
\end{corollary}

\section{Adding a $(2,0)$--large set.}

\begin{theorem}
\label{large}
Assume that $V\subseteq\V^*$ are universes of Set Theory. Let $h\in\baire\cap
\V$ be a strictly increasing function. Suppose that
\[\V^*\models(\exists\eta\in\prod_{n\in\omega}h(n))(\forall\rho\in\prod_{n\in
\omega}\cap\V)(\forall^\infty n\in\omega)(\rho(n)\neq\eta(n)).\]
Then there is a set $X\in [\omega]^{\textstyle \omega}\cap\V^*$ such that
\[\V^*\models(\forall f\in\baire\cap\V)\big(\;(\forall n\in\omega)(n<f(n))\ \
\Rightarrow\ \ \nor \{m\in X: f(m)\in X\}\nor <\omega\big)\]
(so in particular the set $\omega\setminus X$ is $(2,0)$--large over $\V$).
\end{theorem}

\Proof Let $\langle n_i: i\in\omega\rangle$ be defined by
\[n_0=0,\qquad n_{i+1}=n_i+\prod_{k\leq n_i} h(k).\]
Let $H:\bigcup\limits_{i\in\omega}\prod\limits_{k\leq n_i}h(k)\stackrel{1-1}{
\longrightarrow}\omega$ be a bijection such that for each $i\in\omega$
\[H\left[\prod_{k\leq n_i} h(k)\right]= [n_i,n_{i+1}).\]
For a function $f\in\baire$ define $\rho_f\in\prod\limits_{k\in\omega}h(k)$ by
\[\rho_f(k)=\left\{
\begin{array}{ll}
H^{-1}(f(k))(k) &\mbox{if }n_i\leq k<n_{i+1}\mbox{ and }n_{i+1}\leq f(k)\\
0               &\mbox{otherwise.}
\end{array}
\right.\]
Note that the mapping $f\mapsto\rho_f$ is coded in $\V$.\\
Let $X=\{H(\eta\rest n_i): i\in\omega\}$ (so it is an infinite subset of
$\omega$ from $\V^*$). Suppose that $f\in\baire\cap\V$ is such that $(\forall
n\in\omega)(n<f(n))$. Look at $\rho_f$. We know that $\rho_f\in\prod\limits_{
k\in\omega} h(k)\cap\V$. So, by the assumptions on $\eta$, we find
$i_0\in\omega$ such that
\[(\forall i\geq i_0)(\eta(i)\neq\rho_f(i)).\]
Suppose now that $i\geq i_0$ and $f(H(\eta\rest n_i))\in X$. Then $f(H(\eta
\rest n_i))=H(\eta\rest n_j)$ for some $j>i$. But this means that
\[\rho_f(H(\eta\rest n_i))=H^{-1}(H(\eta\rest n_j))(H(\eta\rest n_i))=\eta(
H(\eta\rest n_i)),\]
a contradiction with the choice of $i_0$.~$\QED$

\begin{definition}
Let ${\frak d}(R^\exists_0)$ be the minimal size of a family ${\cal K}$ of
partitions $\langle K_n: n\in\omega\rangle$ of $\omega$ into sets of size
$\geq 2$ such that for every infinite co-infitnite subset $X$ of $\omega$ we
have 
\[(\exists\langle K_n: n\in\omega\rangle\in {\cal K})(\exists^\infty n\in
\omega)(K_n\cap X=\emptyset).\]
\end{definition}
In \cite[3.1]{RoSh:501} we remarked that ${\frak b}\leq {\frak d}(R^\exists_0)
\leq \non(\M)$. Now we may add:

\begin{corollary}
$\kappa^*\leq {\frak d}(R^\exists_0)$.
\end{corollary}

\Proof It follows from \ref{large} (compare the proof of \ref{corcovstar});
remember \ref{thereisno}.~$\QED$

\ifx\undefined\bysame
\newcommand{\bysame}{\leavevmode\hbox to3em{\hrulefill}\,}
\fi

\end{document}